\documentclass[12pt]{article}
\usepackage[]{amsmath,amssymb}
\usepackage{latexsym}
\usepackage{cite}

\newtheorem{definition}{Definition}[section]
\newtheorem{theorem}[definition]{Theorem}
\newtheorem{lemma}[definition]{Lemma}
\newtheorem{corollary}[definition]{Corollary}

\newtheorem{note}[definition]{Note}

\newtheorem{proposition}[definition]{Proposition}

\def\F{\mathbb F}

\begin{document}
\title{\bf How to sharpen a tridiagonal pair\footnote{
Preprint of an article submitted for consideration in Journal of Algebra
and its Applications (JAA)© 2008 copyright World Scientific Publishing Company
http://www.worldscinet.com/jaa/
}}
\author{
Tatsuro Ito\footnote{Supported in part by JSPS grant
18340022.} $\;$   and
Paul Terwilliger\footnote{This author gratefully acknowledges 
support from the FY2007 JSPS Invitation Fellowship Program
for Reseach in Japan (Long-Term), grant L-07512.}
}
\date{}

\maketitle
\begin{abstract}
Let $\F$ denote a
 field and let $V$ denote a vector space over $\F$ with
finite positive dimension.
We consider a pair of linear transformations $A:V \to V$
and $A^*:V \to V$ that satisfy the following conditions:
(i)
each of $A,A^*$ is diagonalizable;
(ii)
there exists an ordering $\lbrace V_i\rbrace_{i=0}^d$ of the eigenspaces of
$A$ such that
$A^* V_i \subseteq V_{i-1} + V_{i} + V_{i+1}$ for $0 \leq i \leq d$,
where $V_{-1}=0$ and $V_{d+1}=0$;
(iii)
there exists an ordering $\lbrace V^*_i\rbrace_{i=0}^\delta$
of the eigenspaces of $A^*$ such that
$A V^*_i \subseteq V^*_{i-1} + V^*_{i} + V^*_{i+1}$ for
 $0 \leq i \leq \delta$,
where $V^*_{-1}=0$ and $V^*_{\delta+1}=0$;
(iv)
there is no subspace $W$ of $V$ such that
$AW \subseteq W$, $A^* W \subseteq W$, $W \neq 0$, $W \neq V$.
We call such a pair a {\it tridiagonal pair} on $V$.
It is known that $d=\delta$, and
for $0 \leq i \leq d$ the dimensions of
$V_i, V^*_i, V_{d-i}, V^*_{d-i}$ coincide.
 Denote
this common dimension by 
$\rho_i$
and call $A,A^*$  {\it sharp} whenever $\rho_0=1$.
Let $T$ denote the $\F$-subalgebra of
$\mbox{\rm End}_\F(V)$ 
generated by $A,A^*$.
We show:
(i) the center $Z(T)$ is a field whose dimension
over $\F$ is $\rho_0$;
(ii) the field $Z(T)$ is isomorphic to
each of
$E_0TE_0$,
$E_dTE_d$, 
$E^*_0TE^*_0$,
$E^*_dTE^*_d$,
where $E_i$ (resp. $E^*_i$) is
the 
 primitive idempotent of $A$ (resp. $A^*$) associated
with $V_i$ (resp. $V^*_i$);
(iii) with respect to the 
$Z(T)$-vector space $V$
the pair $A,A^*$ is a sharp tridiagonal pair.

 \end{abstract}

\section{Tridiagonal pairs}

\noindent 
Throughout this paper $\F$ denotes a field and
$V$ denotes a vector space over $\F$ with finite
positive dimension.

\medskip
\noindent 
We begin by recalling the notion of a tridiagonal pair. 
We will use the following terms.
Let ${\rm End}(V)=
{\rm End}_{\F}(V)$
 denote the $\F$-algebra of all $\F$-linear
transformations from $V$ to $V$.
 Given $A \in {\rm End}(V)$
and a
subspace $W \subseteq V$,
we call $W$ an
 {\it eigenspace} of $A$ whenever 
 $W\not=0$ and there exists $\theta \in \F$ such that 
$W=\lbrace v \in V \vert Av = \theta v\rbrace$;
in this case $\theta$ is the {\it eigenvalue} of
$A$ associated with $W$.
We say that $A$ is {\it diagonalizable} whenever
$V$ is spanned by the eigenspaces of $A$.

\begin{definition}  
{\rm \cite[Definition~1.1]{TD00}}
\label{def:tdp}
\rm
By a {\it tridiagonal pair} (or {\it $TD$ pair})
on $V$
we mean an ordered pair $A,A^*$ of
elements in 
 ${\rm End}(V)$
that satisfy the following four conditions.
\begin{enumerate}
\item Each of $A,A^*$ is diagonalizable.
\item There exists an ordering $\lbrace V_i\rbrace_{i=0}^d$ of the  
eigenspaces of $A$ such that 
\begin{equation}
A^* V_i \subseteq V_{i-1} + V_i+ V_{i+1} \qquad \qquad (0 \leq i \leq d),
\label{eq:t1}
\end{equation}
where $V_{-1} = 0$ and $V_{d+1}= 0$.
\item There exists an ordering $\lbrace V^*_i\rbrace_{i=0}^{\delta}$ of
the  
eigenspaces of $A^*$ such that 
\begin{equation}
A V^*_i \subseteq V^*_{i-1} + V^*_i+ V^*_{i+1} 
\qquad \qquad (0 \leq i \leq \delta),
\label{eq:t2}
\end{equation}
where $V^*_{-1} = 0$ and $V^*_{\delta+1}= 0$.
\item
There does not exist a subspace $W$ of $V$ such  that $AW\subseteq W$,
$A^*W\subseteq W$, $W\not=0$, $W\not=V$.
\end{enumerate}
\end{definition}

\begin{note}
\label{lem:convention}
\rm
According to a common notational convention $A^*$ denotes 
the conjugate-transpose of $A$. We are not using this convention.
In a TD pair $A,A^*$ the linear transformations $A$ and $A^*$
are arbitrary subject to (i)--(iv) above.
\end{note}

\medskip
\noindent We refer the reader to 
\cite{
hasan2,
bas6,
neubauer,
TD00,
shape,
tdanduq,
NN,
qtet,
Ev,
IT:Krawt,
IT:qrac,
nomsharp,
nomtowards,
nomstructure,
nom:mu,
madrid,
Vidar
} and the references therein for background
information on TD pairs.

\medskip
\noindent 
To motivate
 our results we recall some facts about TD pairs.
Let $A,A^*$ denote a TD pair
on $V$, as in Definition 
\ref{def:tdp}. By
\cite[Lemma 4.5]{TD00}
the integers $d$ and $\delta$ from
(ii), (iii) are equal; we call this
common value the {\it diameter} of the
pair.
An ordering of the eigenspaces of $A$ (resp. $A^*$)
is said to be {\em standard} whenever it satisfies 
(\ref{eq:t1})
 (resp. (\ref{eq:t2})). 
We comment on the uniqueness of the standard ordering.
Let $\{V_i\}_{i=0}^d$ denote a standard ordering of the eigenspaces of $A$.
By \cite[Lemma~2.4]{TD00}, 
 the ordering $\{V_{d-i}\}_{i=0}^d$ is also standard and no further
 ordering
is standard.
A similar result holds for the eigenspaces of $A^*$.
Let $\{V_i\}_{i=0}^d$ (resp.
$\{V^*_i\}_{i=0}^d$)
denote a standard ordering of the eigenspaces
 of $A$ (resp. $A^*$).
By \cite[Corollary~5.7]{TD00}, 
for $0 \leq i \leq d$ the spaces $V_i$, $V^*_i$
have the same dimension; we denote
this common dimension by $\rho_i$. 
By \cite[Corollaries 5.7, 6.6]{TD00}
the sequence $\{\rho_i\}_{i=0}^d$ is symmetric and unimodal;
that is $\rho_i=\rho_{d-i}$ for $0 \leq i \leq d$ and
$\rho_{i-1} \leq \rho_i$ for $1 \leq i \leq d/2$.
We call the sequence $\{\rho_i\}_{i=0}^d$ the {\em shape}
of $A,A^*$.
The TD pair $A,A^*$ is called {\it sharp} whenever
$\rho_0=1$.
By
\cite[Theorem~1.3]{nomstructure},
if $\F$ is algebraically closed then
 $A,A^*$ is sharp.
It is an open problem to classify the
sharp TD pairs up to isomorphism,
but progress is being made
\cite{
neubauer,
 IT:Krawt,
NN,
Vidar,
nomtowards,
nom:mu,
IT:qrac
}.
Concerning the nonsharp TD pairs,
relatively little research has been done.
In this paper we get the ball rolling 
by proving three results on the subject.
Referring to the above TD pair $A,A^*$ let 
$T$ denote the $\F$-subalgebra of 
${\rm End}(V)$ generated by $A,A^*$.
We show:
(i) the center $Z(T)$ is a field whose dimension
over $\F$ is $\rho_0$;
(ii) the field $Z(T)$ is isomorphic to
each of
$E_0TE_0$,
$E_dTE_d$, 
$E^*_0TE^*_0$,
$E^*_dTE^*_d$,
where $E_i$ (resp. $E^*_i$) is
the 
 primitive idempotent of $A$ (resp. $A^*$) associated
with $V_i$ (resp. $V^*_i$);
(iii) with respect to the 
$Z(T)$-vector space structure on $V$,
the pair $A,A^*$ is a sharp TD pair.
Our proof is based on the recent discovery
that each of
$E_0TE_0$,
$E_dTE_d$, 
$E^*_0TE^*_0$,
$E^*_dTE^*_d$ is commutative
\cite[Theorem~2.6]{nomstructure}.
Essentially our results follow from
this and the Wedderburn theory
\cite[Chapter~IV]{cr},
but for the sake of completeness
and accessibility
we prove
everything from first principles.

\section{Tridiagonal systems}

\indent
When working with a TD pair, it is often convenient to consider
a closely related object called a TD system.
To define a TD system, we recall a few concepts from linear
algebra.
Let $A$ denote a diagonalizable 
element of $\mbox{\rm End}(V)$.
Let $\{V_i\}_{i=0}^d$ denote an ordering of the eigenspaces of $A$
and let $\{\theta_i\}_{i=0}^d$ denote the corresponding ordering of
the eigenvalues of $A$.
For $0 \leq i \leq d$ define $E_i \in 
\mbox{\rm End}(V)$ 
such that $(E_i-I)V_i=0$ and $E_iV_j=0$ for $j \neq i$ $(0 \leq j \leq d)$.
Here $I$ denotes the identity of $\mbox{\rm End}(V)$.
We call $E_i$ the {\em primitive idempotent} of $A$ corresponding to $V_i$.
Observe that
(i) $I=\sum_{i=0}^d E_i$;
(ii) $E_iE_j=\delta_{i,j}E_i$ $(0 \leq i,j \leq d)$;
(iii) $V_i=E_iV$ $(0 \leq i \leq d)$;
(iv) $A=\sum_{i=0}^d \theta_i E_i$.
Moreover
\begin{eqnarray}         \label{eq:defEi}
  E_i=\prod_{\stackrel{0 \leq j \leq d}{j \neq i}}
          \frac{A-\theta_jI}{\theta_i-\theta_j}.
\end{eqnarray}
Note that each of 
$\{A^i\}_{i=0}^d$,
$\{E_i\}_{i=0}^d$ is a basis for the $\F$-subalgebra
of $\mbox{\rm End}(V)$ generated by $A$.
Moreover $\prod_{i=0}^d(A-\theta_iI)=0$.
Now let $A,A^*$ denote a TD pair on $V$.
An ordering of the primitive idempotents 
 of $A$ (resp. $A^*$)
is said to be {\em standard} whenever
the corresponding ordering of the eigenspaces of $A$ (resp. $A^*$)
is standard.

\begin{definition}
{\rm \cite[Definition~2.1]{TD00}}
 \label{def:TDsystem} 
\rm
By a {\it tridiagonal system} 
(or {\it  $TD$ system}) on $V$ we mean a sequence
 $(A;\{E_i\}_{i=0}^d;A^*;\{E^*_i\}_{i=0}^d)$
that satisfies (i)--(iii) below.
\begin{itemize}
\item[(i)]
$A,A^*$ is a TD pair on $V$.
\item[(ii)]
$\{E_i\}_{i=0}^d$ is a standard ordering
of the primitive idempotents of $A$.
\item[(iii)]
$\{E^*_i\}_{i=0}^d$ is a standard ordering
of the primitive idempotents of $A^*$.
\end{itemize}
\end{definition}

\begin{definition}   
     \label{def}
\rm
Let $\Phi=(A;\{E_i\}_{i=0}^d;A^*$; $\{E^*_i\}_{i=0}^d)$ 
denote a TD system on $V$.
For $0 \leq i \leq d$ let $\theta_i$ (resp. $\theta^*_i$)
denote the eigenvalue of $A$ (resp. $A^*$)
associated with the eigenspace $E_iV$ (resp. $E^*_iV$).
We call $\{\theta_i\}_{i=0}^d$ (resp. $\{\theta^*_i\}_{i=0}^d$)
the {\em eigenvalue sequence}
(resp. {\em dual eigenvalue sequence}) of $\Phi$.
Observe that $\{\theta_i\}_{i=0}^d$ (resp. $\{\theta^*_i\}_{i=0}^d$) are
mutually distinct and contained in $\F$.
By the {\it shape} of $\Phi$ we mean the shape of
the TD pair $A,A^*$.
We say $\Phi$ is
{\it sharp} whenever 
$A,A^*$ is sharp.
\end{definition}
\medskip
\noindent The following characterization of TD systems
will be useful.

\begin{lemma}
\label{lem:char}
A sequence 
$(A; \lbrace E_i\rbrace_{i=0}^d;
A^*;
\lbrace E^*_i\rbrace_{i=0}^d)$
is 
 a TD system on $V$ if and only
if {\rm (i)}--{\rm (iv)} hold below:
\begin{itemize}
\item[\rm (i)]
Each of $A,A^*$ is a diagonalizable element of $\mbox{\rm End}(V)$.
\item[\rm (ii)] $\lbrace E_i\rbrace_{i=0}^d$  is an ordering
of the primitive idempotents of $A$ such that
\begin{eqnarray*}
E_iA^*E_j=0 
\qquad \mbox{\rm if} \quad |i-j|>1, 
\qquad (0 \leq i,j\leq d).
\end{eqnarray*}
\item[\rm (iii)] $\lbrace E^*_i\rbrace_{i=0}^d$  is an ordering
of the primitive idempotents of $A^*$ such that
\begin{eqnarray*}
E^*_iAE^*_j=0 
\qquad \mbox{\rm if} \quad |i-j|>1, 
\qquad (0 \leq i,j\leq d).
\end{eqnarray*}
\item[\rm (iv)] 
There does not exist a subspace $W$ of $V$ such that
$AW\subseteq W$,
$A^*W\subseteq W$,
$W\not=0$,
$W\not=V$.
\end{itemize}
\end{lemma}
\noindent {\it Proof:}  Routine.
\hfill $\Box$ \\

\section{The algebra $T$ and its center}

\noindent For the rest of this paper
fix
a TD system
 $(A;\{E_i\}_{i=0}^d;A^*;\{E^*_i\}_{i=0}^d)$ on $V$,
with eigenvalue sequence
$\lbrace \theta_i\rbrace_{i=0}^d$,
dual eigenvalue sequence
$\lbrace \theta^*_i\rbrace_{i=0}^d$,
and shape $\lbrace \rho_i \rbrace_{i=0}^d$.
Let $T$ denote the $\F$-subalgebra of
 $\mbox{\rm End}(V)$ generated by $A,A^*$.
By definition $T$ contains the identity
$I$ of 
 $\mbox{\rm End}(V)$.
By
(\ref{eq:defEi})
the algebra $T$ contains each of $E_i, E^*_i$
for $0 \leq i \leq d$. By construction
the $T$-module $V$ is faithful.
The $T$-module $V$ is irreducible by
Lemma
\ref{lem:char}(iv).

\begin{definition}
\rm
An element $z \in T$ is called {\it central}
whenever $zt=tz$ for all $t \in T$.
The {\it center} $Z(T)$ is the $\F$-subalgebra
of $T$ consisting of the central elements in $T$.
\end{definition}

\begin{lemma}
$Z(T)$ is a field.
\end{lemma}
\noindent {\it Proof:} 
Since $Z(T)$ is commutative,
it suffices to show 
that each nonzero element $z \in Z(T)$
has an inverse in 
$Z(T)$.
Define $K=\lbrace v \in V\,|\,zv=0\rbrace$
and observe that $K$ is a $T$-submodule of $V$.
Note that $K \not=V$ since $z \not=0$ and
the $T$-module $V$ is faithful.
Now $K=0$ 
since the $T$-module $V$ is 
irreducible.
Therefore $z$ is invertible on $V$.
Let $z^{-1} \in {\rm End}(V)$
denote the inverse of $z$ on $V$.
By elementary linear algebra
$z^{-1}$ is a polynomial in $z$
and is therefore contained in $Z(T)$.
The result follows.
\hfill $\Box$ \\

\begin{definition}
\label{def:rhodim}
\rm The center $Z(T)$ is an
$\F$-subalgebra of $T$ and therefore
an $\F$-vector space.
Let $\rho$ denote the dimension of the
$\F$-vector space $Z(T)$. 
\end{definition}

\noindent In Lemma
\ref{lem:last}
 we will show
that 
the parameter
 $\rho$
 from Definition
\ref{def:rhodim} is equal to
$\rho_0$.

\medskip
\noindent For the $T$-module $V$
the restriction of the $T$-action
to $Z(T)$ gives a $Z(T)$-vector space
structure on $V$. 
Each of $A,A^*$ 
commutes with
everything in $Z(T)$ and 
is therefore a $Z(T)$-linear transformation.
Let 
 ${\mbox{\rm End}}_{Z(T)}(V)$ denote the
$Z(T)$-algebra consisting of the
$Z(T)$-linear transformations from $V$ to $V$. 
By the construction
$A,A^* \in 
 {\rm End}_{Z(T)}(V)$.

\medskip
\noindent Of course every $Z(T)$-subspace of $V$
is an $\F$-subspace of $V$,
but  not every $\F$-subspace of
 $V$
is a $Z(T)$-subspace of $V$.

\begin{lemma}
\label{lem:goodsubdim}
For each $Z(T)$-subspace $W$ of $V$,
\begin{eqnarray*}
\rho \, \mbox{\rm dim}_{Z(T)}(W) = \mbox{\rm dim}_\F(W).
\end{eqnarray*}
\end{lemma}
\noindent {\it Proof:} 
$W$ is the direct sum of
$\mbox{\rm dim}_{Z(T)}(W)$ many
one-dimensional $Z(T)$-subspaces of $V$. By
Definition
\ref{def:rhodim},
each of these
has dimension $\rho$
as an $\F$-subspace of $V$.
The result follows.
\hfill $\Box$ \\

\begin{lemma}
\label{lem:diag}
Let $\psi$ denote any element of 
${\rm End}_{Z(T)}(V)$.
If $\psi$
is diagonalizable when viewed as an $\F$-linear transformation
on the $\F$-vector space $V$, then
 $\psi$ is diagonalizable when viewed as a 
$Z(T)$-linear transformation on the 
$Z(T)$-vector space $V$. In this case
the eigenspaces, eigenvalues, and
primitive idempotents of $\psi$ are 
the same for the two points of view.
\end{lemma}
\noindent {\it Proof:} 
For $\theta \in \F$ the 
set $\lbrace v \in V \,|\,\psi v = \theta v\rbrace$
is the same for the two points of view.
\hfill $\Box$ \\

\begin{lemma}
\label{lem:irred2}
There does not exist a $Z(T)$-subspace $W$ of
$V$ such that
$AW \subseteq W$,
$A^*W \subseteq W$,
$W\not=0$,
$W \not=V$.
\end{lemma}
\noindent {\it Proof:} 
By
Lemma
\ref{lem:char}(iv) there does not
exist an $\F$-subspace $W$ of $V$
such that
$AW \subseteq W$,
$A^*W \subseteq W$,
$W\not=0$,
$W \not=V$.
The result follows since each $Z(T)$-subspace of
$V$ is an $\F$-subspace of $V$.
\hfill $\Box$ \\

\begin{proposition}
\label{thm:tdk}
The elements $(A;\lbrace E_i\rbrace_{i=0}^d;
A^*;\lbrace E^*_i\rbrace_{i=0}^d)$
act on the $Z(T)$-vector space
$V$ as a TD system.
This TD system has eigenvalue sequence 
$\lbrace \theta_i\rbrace_{i=0}^d$,
dual eigenvalue sequence
$\lbrace \theta^*_i\rbrace_{i=0}^d$,
and shape $\lbrace \rho_i/\rho\rbrace_{i=0}^d$.
\end{proposition}
\noindent {\it Proof:} 
To verify the first assertion, we check that
$(A;\lbrace E_i\rbrace_{i=0}^d;
A^*;\lbrace E^*_i\rbrace_{i=0}^d)$
satisfy the conditions
of
Lemma
\ref{lem:char}.
Concerning condition (i),
we mentioned 
below Definition
\ref{def:rhodim}
that  
each of $A,A^*$ is
contained in
$\mbox{\rm End}_{Z(T)}(V)$.
Using Lemma  
\ref{lem:diag} we find that
the map $A$ (resp. $A^*$) is diagonalizable
as a 
$Z(T)$-linear transformation
on the $Z(T)$-vector space $V$,
 since
$A$ (resp. $A^*$) is 
 diagonalizable
as an $\F$-linear transformation on the
$\F$-vector space $V$.
Using Lemma \ref{lem:diag}
we find that conditions (ii), (iii) hold
when $A$, $A^*$ are viewed as 
$Z(T)$-linear transformations on the
$Z(T)$-vector space $V$,
since they hold
when $A,A^*$ are viewed as
$\F$-linear transformations on the
  $\F$-vector space $V$.
Condition (iv) 
holds by Lemma
\ref{lem:irred2}.
Our first assertion is now verified.
The assertions
about the eigenvalue sequence
and dual eigenvalue sequence
follow from
Definition
     \ref{def}
and Lemma
\ref{lem:diag}.
The assertion about the shape follows from
Lemma
\ref{lem:goodsubdim} and since for
$0 \leq i \leq d$,
each of $E_iV$, $E^*_iV$ is a $Z(T)$-subspace of
$V$ which 
has dimension
$\rho_i$ as an $\F$-subspace of $V$.
\hfill $\Box$ \\

\section{The algebras 
$E_0TE_0$,
$E_dTE_d$,
$E^*_0TE^*_0$,
$E^*_dTE^*_d$}

\noindent
In this section we show how
the center $Z(T)$ is related to 
$E_0TE_0$,
$E_dTE_d$,
$E^*_0TE^*_0$,
$E^*_dTE^*_d$. 
Invoking symmetry we will restrict our attention to
$E^*_0TE^*_0$.
We view $E^*_0TE^*_0$ as an $\F$-algebra with
multiplicative identity $E^*_0$. 
Our first goal is to 
show that this algebra is a field.
Our point of departure is the following recent discovery.

\begin{lemma}
\label{lem:comgen}
{\rm \cite[Theorem~2.6]{nomstructure}}
The algebra $E^*_0TE^*_0$ is commutative and generated by
\begin{eqnarray*}
  E^*_0A^iE^*_0 \qquad \qquad (1\leq i \leq d).
\end{eqnarray*}
\end{lemma}

\begin{lemma}
\label{lem:ete0}
The $E^*_0TE^*_0$-module $E^*_0V$ is
irreducible.
\end{lemma}
\noindent {\it Proof:} 
Let $W$ denote a nonzero 
$E^*_0TE^*_0$-submodule of 
$E^*_0V$.
 We show
that
$W=E^*_0V$. 
Pick any nonzero $w \in W$ and
observe
$E^*_0TE^*_0w \subseteq  W$.
The element $E^*_0$ acts on
$E^*_0V$ as the identity so
$E^*_0Tw \subseteq W$.
Observe that $Tw$ is a 
nonzero $T$-submodule of $V$.
The $T$-module $V$ is irreducible 
so 
$Tw=V$.
Therefore 
$E^*_0Tw=E^*_0V$.
By the above comments
$W=E^*_0V$ and the result follows.
\hfill $\Box$ \\

\begin{lemma}
\label{lem:es0inj}
For all nonzero $s \in TE^*_0$
the restriction of $s$ to
$E^*_0V$ is injective.
\end{lemma}
\noindent {\it Proof:} 
We define
$K=\lbrace v \in E^*_0V |sv=0\rbrace$
and show $K=0$. In view of
Lemma
\ref{lem:ete0}
it suffices to show that $K$ is invariant
under $E^*_0TE^*_0$, and that
$K \not=E^*_0V$.
We now show that $K$ is invariant under
$E^*_0TE^*_0$. 
For all $k \in K$ and $t \in T$ we show
$E^*_0tE^*_0k \in K$.
By construction
$E^*_0tE^*_0 k \in E^*_0V$,
so we just have to check that
$s E^*_0tE^*_0k = 0$.
Suppose $s E^*_0tE^*_0k \not= 0$.
 Then $Ts E^*_0tE^*_0k = V$
by the irreducibility of the 
$T$-module $V$, so
$E^*_0Ts E^*_0tE^*_0k = E^*_0V$.
Since $E^*_0TE^*_0$ is commutative and
$Ts \subseteq T$,
each element of
$E^*_0Ts E^*_0$
commutes with 
$E^*_0tE^*_0$.
Using this and $E^{*2}_0 = E^*_0$
we obtain
\begin{eqnarray}
\label{eq:flip}
E^*_0tE^*_0TsE^*_0k = E^*_0V.
\end{eqnarray}
We have  $sE^*_0=s$ since $s \in TE^*_0$, 
and $sk=0$ by construction, so the left-hand
side of
(\ref{eq:flip}) is zero.
Of course the right-hand side of
(\ref{eq:flip}) is nonzero
and we have a contradiction. Therefore
$sE^*_0tE^*_0k=0$, so
$E^*_0tE^*_0k \in K$. 
We have now shown that $K$ is invariant under
$E^*_0TE^*_0$.
We now show
that
$K \not=E^*_0V$.
Suppose $K=E^*_0V$.
Then $s$ is zero on
$E^*_0V$ so $sE^*_0=0$.
But $sE^*_0=s$, so
$s=0$ 
 for  
a contradiction. Therefore $K\not=E^*_0V$.
We have shown
 that $K$ is invariant
under $E^*_0TE^*_0$, and that
$K \not=E^*_0V$.
Therefore $K=0$ in view of 
Lemma
\ref{lem:ete0}.
The result follows.
\hfill $\Box$ \\

\medskip
\noindent Observe that $TE^*_0$ has a $T$-module
structure such that $t.s=ts$ for all $t \in T$
and $s \in TE^*_0$.

\begin{lemma}
\label{lem:vbij}
For all nonzero $v \in E^*_0V$ the map
\begin{eqnarray*}
 TE^*_0  & \to &   V \\
 s  \;\; &\mapsto &   sv
\label{eq:ssv}
\end{eqnarray*}
is an isomorphism of $T$-modules.
\end{lemma}
\noindent {\it Proof:} 
By construction the map
is a homomorphism
of $T$-modules.
The map is injective by
Lemma
\ref{lem:es0inj}.
To see that the map is surjective,
note that its image is a $T$-submodule
of $V$, and nonzero since it contains $v$.
This image is equal to  
$V$ by the irreducibility of
the $T$-module $V$.
\hfill $\Box$ \\

\noindent We emphasize one point
from Lemma
\ref{lem:vbij}.

\begin{corollary}
\label{cor:tes}
The $T$-module $TE^*_0$ is irreducible and
faithful.
\end{corollary}
\noindent {\it Proof:} 
Immediate from Lemma
\ref{lem:vbij} and since the
$T$-module $V$ is irreducible and faithful.
\hfill $\Box$ \\

\begin{proposition}
\label{prop:ext}
$E^*_0TE^*_0$ is a field.
\end{proposition}
\noindent {\it Proof:} 
Since $E^*_0TE^*_0$ is commutative,
it suffices to show that each nonzero 
   $a \in 
 E^*_0TE^*_0$
has an inverse in 
 $E^*_0TE^*_0$.
Since $a$ is nonzero and contained in 
$TE^*_0$, the space $Ta$ is a nonzero $T$-submodule
of $TE^*_0$. The $T$-module
$TE^*_0$ is irreducible so
$Ta=TE^*_0$. Therefore there exists
$t \in T$ such that
$ta=E^*_0$. Define $b=E^*_0tE^*_0$ and observe $b \in E^*_0TE^*_0$.
Since $a \in E^*_0TE^*_0$ we find $E^*_0a=a$;
using this and $ta=E^*_0$ 
 we find
$ba=E^*_0$. Therefore $b$ is the desired inverse of
$a$ and 
the result follows.
\hfill $\Box$ \\

\noindent Our next goal is to show that the
fields $Z(T)$, $E^*_0TE^*_0$ are isomorphic.
To this end it is helpful to define a certain bilinear
form.

\begin{definition}
\label{def:bilform}
\rm
Observe that $TE^*_0$ (resp. $E^*_0T$) has a right (resp. left)
$E^*_0TE^*_0$-module
structure.
Using these structures we view each of
$TE^*_0, E^*_0T$ as a vector space over
the field $E^*_0TE^*_0$.
Over this field we define a bilinear form
$(\,,\,): E^*_0T \times TE^*_0 \rightarrow E^*_0TE^*_0$
such that $(r, s)=rs$ 
for all $r \in E^*_0T$ and $s \in TE^*_0$.
\end{definition}

\begin{lemma}
\label{lem:ndeg}
The bilinear form 
$(\,,\,)$ 
from Definition \ref{def:bilform}
is nondegenerate.
\end{lemma}
\noindent {\it Proof:} 
Given $r \in E^*_0T$ 
such that 
$(r,s)=0$ for all $s \in TE^*_0$,
we show $r=0$. 
By the construction $rTE^*_0=0$. Now $r=0$ since
the action of $T$ on $TE^*_0$ is faithful.
Given $s \in TE^*_0$ such that
$(r,s)=0$ for all $r \in E^*_0T$,
we show $s=0$.
By construction
$E^*_0Ts=0$.
 Suppose $s \not=0$.
Then $Ts$ is a nonzero $T$-submodule
of $TE^*_0$, so $Ts=TE^*_0$ by
the irreducibility of the $T$-module
$TE^*_0$.
Therefore 
$E^*_0Ts=E^*_0TE^*_0$. In this equation
the left hand side is zero and the right
hand side is nonzero, for a contradiction.
Hence $s=0$
and 
the result follows.
\hfill $\Box$ \\

\begin{definition}
\label{def:dualb}
\rm
By Lemma
\ref{lem:ndeg} the
 $E^*_0TE^*_0$-vector spaces
  $TE^*_0$ and $E^*_0T$  have the same dimension
which we denote by $n$.
Let $\lbrace x_i \rbrace_{i=1}^n$ denote
a basis for the 
 $E^*_0TE^*_0$-vector space
$TE^*_0$ and let
$\lbrace x'_i \rbrace_{i=1}^n$ denote
the basis for
the  
 $E^*_0TE^*_0$-vector space
$E^*_0T$ which is dual to $\lbrace x_i\rbrace_{i=1}^n$
with respect to
$(\,,\,)$.
By construction
\begin{eqnarray}
\label{eq:bildual}
 x'_ix_j=\delta_{i,j}E^*_0
\qquad \qquad (1 \leq i,j\leq n).
\end{eqnarray}
\end{definition}

\begin{lemma}
\label{lem:centact}
The map
\begin{equation}
\begin{split}
 Z(T)\; \; & \to  \;\;  E^*_0TE^*_0 \\
 z\; \; \;\; &\mapsto  \;\; zE^*_0 
\label{eq:cent}
\end{split}
\end{equation}
is an isomorphism of fields.
\end{lemma}
\noindent {\it Proof:} 
For $z \in Z(T)$ we have 
$zE^*_0=E^*_0zE^*_0$, since
$E^{*2}_0=E^*_0$ and
$zE^*_0=E^*_0z$.
Therefore
$zE^*_0 \in E^*_0TE^*_0$.
By construction the map 
defined in line (\ref{eq:cent})
is an $\F$-algebra homomorphism.
We now show that this map is a bijection.
The map 
is injective 
since a field has no nonzero proper ideals.
To see that the map 
is surjective,
for a given $a \in E^*_0TE^*_0$
we display $z \in Z(T)$ such that $zE^*_0=a$.
Define $z=\sum_{i=1}^n x_ia x'_i$ where
$\lbrace x_i\rbrace_{i=1}^n$ and $\lbrace x'_i\rbrace_{i=1}^n$ 
are from
Definition
\ref{def:dualb}.
Using
(\ref{eq:bildual}) and $a=aE^*_0$ we obtain
$zx_j =x_ja$ for $1 \leq j \leq n$.
Recall $\lbrace x_j\rbrace_{j=1}^n$ is a basis for
the $E^*_0TE^*_0$-vector space $TE^*_0$,
so $zs=sa$ for all $s \in TE^*_0$. 
In other words $z$ acts on $TE^*_0$ as 
$s \to sa$.
Using this we find that for all $t \in T$ the expression
$zt-tz$ vanishes on $TE^*_0$,
yielding $zt=tz$ since the $T$-module $TE^*_0$ is faithful.
Therefore $z \in Z(T)$.
We mentioned above that $zs=sa$ for all $s \in TE^*_0$;
taking $s=E^*_0$ we find
$zE^*_0=E^*_0a=a$.
We have shown $z \in Z(T)$ and
$zE^*_0=a$, 
so the
the map
is surjective.
 The result follows.
\hfill $\Box$ \\

\noindent We need one fact about dimensions and then we will be 
ready for the main result.

\begin{lemma} 
\label{lem:last}
With reference to Definition
\ref{def:rhodim}, we have
$\rho=\rho_0$.
\end{lemma}
\noindent {\it Proof:} 
By Definition
\ref{def:rhodim} the parameter  $\rho$ is the
dimension of the $\F$-vector space
$Z(T)$.
The map $Z(T) \to E^*_0TE^*_0$ from
Lemma
\ref{lem:centact}
is an isomorphism of $\F$-vector spaces,
so 
the $\F$-vector spaces 
$Z(T)$, $E^*_0TE^*_0$ have the same dimension.
The map
$TE^*_0\to V$
from Lemma
\ref{lem:vbij}
is an isomorphism of $T$-modules
and hence an isomorphism of 
$\F$-vector spaces.
Under this map
the image of $E^*_0TE^*_0$ is $E^*_0V$,
so the $\F$-vector spaces $E^*_0TE^*_0$ and $E^*_0V$ have the
same dimension.
By construction the $\F$-vector space $E^*_0V$ has dimension $\rho_0$.
The result follows.
\hfill $\Box$ \\

\noindent Combining our above results we immediately obtain
the following.

\begin{theorem}
\label{thm:final}
 Let $(A;\{E_i\}_{i=0}^d;A^*;\{E^*_i\}_{i=0}^d)$ 
denote a TD system 
on the $\F$-vector space $V$,
with eigenvalue sequence
$\lbrace \theta_i\rbrace_{i=0}^d$,
dual eigenvalue sequence
$\lbrace \theta^*_i\rbrace_{i=0}^d$,
and shape $\lbrace \rho_i \rbrace_{i=0}^d$.
Let $T$ denote the $\F$-subalgebra of
 $\mbox{\rm End}_{\F}(V)$ generated by $A,A^*$. 
Then the following
{\rm (i)}--{\rm (iii)} hold.
\begin{itemize}
\item[\rm (i)]
The center $Z(T)$ is a field whose dimension
over $\F$ is $\rho_0$.
\item[\rm (ii)]
The field $Z(T)$ is isomorphic to
each of
$E_0TE_0$,
$E_dTE_d$, 
$E^*_0TE^*_0$,
$E^*_dTE^*_d$.
\item[\rm (iii)]
The elements $(A;\lbrace E_i\rbrace_{i=0}^d;
A^*;\lbrace E^*_i\rbrace_{i=0}^d)$
act on the $Z(T)$-vector space
$V$ as a TD system  that has
 eigenvalue sequence 
$\lbrace \theta_i\rbrace_{i=0}^d$,
dual eigenvalue sequence
$\lbrace \theta^*_i\rbrace_{i=0}^d$,
and shape $\lbrace \rho_i/\rho_0\rbrace_{i=0}^d$.
In particular this TD system is sharp.
\end{itemize}
\end{theorem}

\section{Acknowledgements}
The authors thank
Eric Egge,
Darren Funk-Neubauer,
Kazumasa Nomura, and
Arlene Pascasio
for giving this paper a close reading and offering many
valuable suggestions.

\bigskip


\noindent Tatsuro Ito \hfil\break
\noindent Division of Mathematical and Physical Sciences \hfil\break
\noindent Graduate School of Natural Science and Technology\hfil\break
\noindent Kanazawa University \hfil\break
\noindent Kakuma-machi,  Kanazawa 920-1192, Japan \hfil\break
\noindent email:  {\tt tatsuro@kenroku.kanazawa-u.ac.jp}

\bigskip

\noindent Paul Terwilliger \hfil\break
\noindent Department of Mathematics \hfil\break
\noindent University of Wisconsin \hfil\break
\noindent 480 Lincoln Drive \hfil\break
\noindent Madison, WI 53706-1388 USA \hfil\break
\noindent email: {\tt terwilli@math.wisc.edu }\hfil\break


\begin{thebibliography}{10}

\bibitem{hasan2}
H.~Alnajjar  and B.~Curtin.
\newblock
A family of tridiagonal pairs related to
the quantum affine algebra 
$U\sb q(\widehat{\mathfrak{sl}}\sb 2)$.
\newblock {\em
Electron. J. Linear Algebra}
{\bf 13} (2005) 1--9. 


\bibitem{bas6}
P.~Baseilhac.
\newblock
A family of tridiagonal pairs and related symmetric functions. 
\newblock
{\em J. Phys. A} {\bf 39}  (2006) 11773--11791. 


\bibitem{cr}
 C.~Curtis and I. Reiner.
\newblock
{\em  Representation theory of finite groups and associative algebras}, 
\newblock
Pure and Applied Mathematics, Vol. XI Interscience Publishers,
 a division of John Wiley \& Sons, New York-London 1962.


\bibitem{neubauer}
D. Funk-Neubauer.
\newblock 
   Tridiagonal pairs and the $q$-tetrahedron algebra.
\newblock{\em Linear Algebra Appl.}, submitted for publication;
{\tt arXiv:0806.0901}.




\bibitem{TD00}
T.~Ito, K.~Tanabe, and P.~Terwilliger.
\newblock Some algebra related to ${P}$- and ${Q}$-polynomial association
  schemes,  in:
\newblock {\em Codes and Association Schemes (Piscataway NJ, 1999)}, Amer.
Math. Soc., Providence RI, 2001, pp.
     167--192; 
{\tt arXiv:math.CO/0406556}.

\bibitem{shape}
T.~Ito and P.~Terwilliger.
\newblock The shape of a tridiagonal pair.
\newblock {\em J. Pure Appl. Algebra}
	      {\bf 188}
	    (2004)
		     145--160;
{\tt arXiv:math.QA/0304244}.

\bibitem{tdanduq}
T.~Ito and P.~Terwilliger.
\newblock {Tridiagonal pairs and the quantum affine 
algebra
$U_q({\widehat{sl}}_2)$.}
\newblock {\em Ramanujan J.} 
{\bf 13} (2007) 39--62;
{\tt arXiv:math.QA/0310042}.

\bibitem{NN}
T.~Ito and P.~Terwilliger.
\newblock Two non-nilpotent linear transformations that
satisfy the cubic $q$-Serre relations.
\newblock{\em 
J. Algebra Appl.} {\bf 6} (2007) 477--503;
{\tt arXiv:math.QA/0508398}.

\bibitem{qtet}
T.~Ito and P.~Terwilliger.
\newblock The $q$-tetrahedron algebra and its
finite-dimensional irreducible modules.
\newblock{\em
 Comm. Algebra}  {\bf 35}  (2007) 3415--3439;
{\tt arXiv:math.QA/0602199}.


\bibitem{Ev}
T.~Ito and P.~Terwilliger.
\newblock Finite-dimensional irreducible modules for
the three-point $\mathfrak{sl}_2$ loop algebra.
\newblock{\em 
 Comm. Algebra},
  in press;
{\tt arXiv:0707.2313}. 

\bibitem{IT:Krawt}
T.~Ito and P.~Terwilliger.
\newblock{Tridiagonal pairs of Krawtchouk type.}
\newblock{\em
Linear Algebra Appl.} {\bf 427}  (2007) 218--233;
{\tt arXiv:0706.1065}.


\bibitem{IT:qrac}
T.~Ito and P.~Terwilliger.
\newblock
Tridiagonal pairs of $q$-Racah type.
\newblock{\em
J. Algebra},
submitted for publication;
{\tt arXiv:0807.0271}.



\bibitem{nomsharp}
K.~Nomura and P.~Terwilliger.
\newblock
Sharp tridiagonal pairs.
\newblock {\em Linear Algebra Appl.},
{\bf 429} (2008) 79--99;
{\tt arXiv:0712.3665}. 


\bibitem{nomtowards}
K.~Nomura and P.~Terwilliger.
\newblock
Towards a classification of the tridiagonal pairs.
\newblock {\em Linear Algebra Appl.},
{\bf 429} (2008) 503--518;
{\tt arXiv:0801.0621}. 


\bibitem{nomstructure}
K.~Nomura and P.~Terwilliger.
\newblock
The structure of a tridiagonal pair.
\newblock {\em Linear Algebra Appl.},
in press;
{\tt arXiv:0802.1096}. 


\bibitem{nom:mu}
K.~Nomura and P.~Terwilliger.
\newblock
Tridiagonal pairs and the $\mu$-conjecture.
\newblock {\em Linear Algebra Appl.},
submitted for publication.





\bibitem{madrid}
P.~Terwilliger.
\newblock
An algebraic approach to the Askey scheme of orthogonal polynomials. 
Orthogonal polynomials and special functions, 
255--330, Lecture Notes in Math., 1883, 
Springer, Berlin, 2006; 
{\tt arXiv:math.QA/0408390}. 


\bibitem{Vidar}
M.~Vidar,
\newblock 
Tridiagonal pairs of shape $(1,2,1)$.
\newblock {\em Linear Algebra Appl.}, 
{\bf 429} (2008) 403-428;
{\tt arXiv:0802.3165}. 

 \end{thebibliography}
\end{document}